\documentclass{amsart}
\usepackage[utf8]{inputenc}
\usepackage{amsmath}
\usepackage{amssymb}
\usepackage{amsthm}
\usepackage{color}
\usepackage{graphicx}
\usepackage{subeqnarray}
\usepackage[active]{srcltx}

\newtheorem{theorem}{\bf Theorem}[section]

\newtheorem{rem}{\bf Remark}[section]
\newcommand{\R}{\mathbb{R}}
\makeatletter \renewcommand\@biblabel[1]{#1.} \makeatother

\begin{document}

\title{A note on the steady Poiseuille flow of Carreau-Yasuda fluid}
\author[N. Kutev]{N. Kutev}
\address{Nikolay Kutev, Institute of Mathematics and Informatics, Bulgarian Academy of Sciences, Acad. G. Bontchev str., bl. 9, 1113 Sofia, Bulgaria}
\email{kutev@math.bas.bg}

\author[S. Tabakova]{S. Tabakova*}
\address{Sonia Tabakova, Institute of Mechanics, Bulgarian Academy of Sciences, Acad. G. Bontchev str., bl. 4, 1113 Sofia, Bulgaria}
\email{stabakova@gmail.com}
\thanks{*corresponding author}



\date{}

\begin{abstract}
The steady Poiseuille flow of Carreau-Yasuda fluid in a pipe, caused by constant pressure gradient, is studied theoretically. It is proved that at some values of the viscosity model parameters, the problem has a classical solution, while at others - generalized solution. For the latter, a necessary and sufficient condition is found, which depends on the pressure gradient and Carreau-Yasuda model parameters.
\end{abstract}

\maketitle

\section{Introduction}

The Poiseuille flow problem is one of the fundamental problems of fluid mechanics, corresponding to a parabolic velocity profile in a pipe or channel, when the flow is laminar and fully developed in axial direction. For Newtonian fluid, the problem has a well known analytical solution \cite{Schlichting:2000}. However, for non-Newtonian fluids, it has no analitical solution, except for some special cases for generalized Newtonian fluid models, such as power law model of viscosity dependence on shear rate \cite{Bird:2002}. The behavior of the so-called shear-thinning fluids (viscosity is decreasing function of shear rate), for example polymer solutions, polymer melts, suspensions, emulsions, and some biological fluids, is described by different viscosity models, such as the power law model, Carreau model, Carreau-Yasuda model and others \cite{Bird:2002}-\cite{Boyd:2007}.
 
The general unsteady and oscillatory cases of the Poiseuille flow in an infinite channel or pipe using the Carreau or Carreau-Yasuda models are studied in our previous papers 
\cite{Tabakova:2020} - \cite{Kutev:2021c}.  
The most interesting case of negative power index in these models, $n<0$, corresponds to severe pseudoplastic behaviour \cite{Kotsilkova:2021}-\cite{Brabazon:2003}. It  has been studied theoretically in \cite{Kutev:2021b}, \cite{Kutev:2021c}, but the sufficient and necessary condition for solution existence is still a challenge. The aim of the present work is to prove the solution existence or non-existence for the different parameters of the steady Poiseuille flow problem of Carreau-Yasuda fluid and to find the existence condition at $n<0$.
 
The dimensionless equations of motion and continuity for the steady pipe flow are:
\begin{equation}
\textbf{v}\cdot\nabla \textbf{v}=-\nabla p+\nabla\cdot \textbf{T},
\label{eq_1}
\end{equation}
\begin{equation}
\nabla\cdot\textbf{v}=0,
\label{eq_2}
\end{equation}
where $\textbf{v}$ is velocity vector, $p$ - pressure, $\textbf{T}=\mu_{app}\dot{S}$ - viscous stress tensor with $\displaystyle{\dot{S}=\frac{1}{2}\left(\nabla \textbf{v}+\nabla \textbf{v}^T \right)}$ as strain rate tensor and $\mu_{app}$ - apparent viscosity (constant for Newtonian fluid and a non-linear function of shear rate $\displaystyle{\dot{\gamma}=\sqrt{2\dot{S}:\dot{S}}}$ for generalized Newtonian fluid). 

The flow is assumed fully developed laminar driven by a constant gradient $b$ in axial direction. Then, in the adopted cylindrical coordinate system $(x, Y, \varphi)$ with $x$ as the axial coordinate, the flow velocity is simplified as $\textbf{v}=(U(Y), 0,0)$ and $\textbf{T}$ reduces to the scalar $\mu_{app}(\dot{\gamma})\dot{\gamma}$, where $\displaystyle{\dot{\gamma}=\mid U_Y \mid}$. 
The Carreau-Yasuda viscosity model is assumed for the apparent viscosity
  $\mu_{app}$ :
\begin{equation}
\mu_{app}=1-c+c[1+Cu^\alpha\dot{\gamma}^\alpha]^{(n-1)/\alpha},
\label{eq_4}
\end{equation}
where $b$ is the dimensionless pressure gradient and $\displaystyle{c =1-\frac{\mu_\infty}{\mu_{0}}}$, with $\mu_0$, $\mu_\infty$ as the upper and lower limits of the viscosity corresponding to the low and high shear rates, $Cu$ is the Carreau number (Weissenberg number), $\alpha$ and $n$ are empirically determined for the considered fluid.

Applying the upper assumptions, the system (\ref{eq_1}) - (\ref{eq_4}) is reduced to one elliptic equation for the axial velocity:
\begin{equation}\label{1.1}
L(U(Y))=\frac{1}{Y}\frac{d}{dY}\Bigg\lbrace\left [1-c+c \left(1+Cu^\alpha \mid U_Y \mid ^\alpha\right)^{\frac{n-1}{\alpha}}\right]Y U_{Y}\Bigg\rbrace=b
\end{equation}
with the axisymmetrical and no-slip boundary conditions: 
\begin{equation}\label{1.2}
U_Y(0) = U(R) = 0.  
\end{equation}
where $R$ is the dimensionless pipe radius (the pipe radius or diameter are usually taken as characteristic length).
Here $Cu>0$, $b\in \R$, $\alpha>0$, $n\in \mathbb{R}$, $c\in[0,1]$ are empirically determined parameters. 

If $b=0$, then (\ref{1.1}), (\ref{1.2}) has only the trivial solution $U(Y)\equiv 0$. It is clear that $U(Y,-b)=-U(Y,b)$, where $U(Y,b)$, $U(Y,-b)$ are solutions of (\ref{1.1}), (\ref{1.2}) with right hand side $b$ and $-b$, correspondingly. That is why, further on we consider only the case $b>0$.

The outline of the paper is as follows. In section 2 the main results are formulated in four theorems to be proved in section 3. Two applications, presented by other three theorems,  together with some illustrations of the theorems' results are found in section 4. Finally, comparison between the criteria proved for the unsteady and steady case are given as an appendix. 

\section{Main results}
In order to formulate the main results, we introduce the function
\begin{equation}\label{2.1}
F(\zeta)=\left[1-c+c\left(1+Cu^\alpha\zeta^\alpha\right)^{\frac{n-1}{\alpha}}\right]\zeta \quad \textrm{for} \quad \zeta \geq 0,
\end{equation}
where $\displaystyle{\zeta=\mid U_Y \mid}$.
\begin{theorem}\label{T2.1}
Suppose $b>0$, $\alpha>0$ and either $c=0$, $Cu=0$ or $n=1$. Then problem \eqref{1.1}, \eqref{1.2} is Newtonian one (classical Poiseuille problem in a pipe) and has a unique classical solution:
\begin{equation}\label{2.2}
U(Y)=\frac{b}{4}(Y^2-R^2) \quad \textrm{for} \quad every \quad Y\in [0,R].
\end{equation}
\end{theorem}
\begin{theorem}\label{T2.2}
Suppose $b>0$, $\alpha>0$, $c\in(0,1)$ and either one of the following conditions holds: 
\begin{equation}\label{2.3}
n\geq 0;
\end{equation}
\begin{equation}\label{2.4}
n<0 \quad and \quad \left(1-\frac{\alpha+1}{n}\right)^{\frac{n-1-\alpha}{\alpha}}<\frac{1-c}{\alpha c}.
\end{equation}
Then problem \eqref{1.1}, \eqref{1.2} has a unique classical solution $U(Y)\in C^2([0,R])$ satisfying
\begin{equation}\label{2.5}
U(Y)=-\int^R_Y F^{-1}\left(\frac{bs}{2}\right)ds \quad \textrm{and} 
\end{equation}
\begin{equation}\label{2.6}
0\leq U_Y(Y)\leq\ F^{-1}\left(\frac{bR}{2}\right) \quad \textrm{for} \quad Y\in [0,R].
\end{equation}
\end{theorem}
\begin{theorem}\label{T2.3}
Suppose $b>0$, $\alpha>0$, $n<0$, $c\in(0,1)$ and 
\begin{equation}\label{2.6a}
\left(1-\frac{\alpha+1}{n}\right)^{\frac{n-1-\alpha}{\alpha}}=\frac{1-c}{\alpha c}
\end{equation}
If
\begin{equation}\label{2.7}
\left[1-c+c\left(1-\frac{\alpha+1}{n}\right)^{\frac{n-1}{\alpha}}\right]Cu^{-1}\left(-\frac{\alpha+1}{n}\right)^{\frac{1}{\alpha}}\geq \frac{bR}{2},
\end{equation}
then problem \eqref{1.1}, \eqref{1.2} has a unique classical solution $U(Y)$ satisfying \eqref{2.5}, \eqref{2.6};\\
If
\begin{equation}\label{2.8}
\left[1-c+c\left(1-\frac{\alpha+1}{n}\right)^{\frac{n-1}{\alpha}}\right]Cu^{-1}\left(-\frac{\alpha+1}{n}\right)^{\frac{1}{\alpha}}< \frac{bR}{2},
\end{equation}
then problem \eqref{1.1}, \eqref{1.2} has a unique generalized solution $U(Y)\in C^2([0,R] \setminus {Y_0})\bigcap C^1([0,R])$ satisfying \eqref{2.5}, \eqref{2.6}, $U_{YY}(Y_0)=\infty$, where
\begin{equation}\label{2.9}
Y_0=\frac{2}{b}\left[1-c+c\left(1-\frac{\alpha+1}{n}\right)^{\frac{n-1}{\alpha}}\right]Cu^{-1}\left(-\frac{\alpha+1}{n}\right)^{\frac{1}{\alpha}} \quad \in(0,R).
\end{equation}
\end{theorem}
\begin{rem}\label{R2.1}
If \eqref{2.7} is strict inequality then $U(Y)\in C^2([0,R])$, while when \eqref{2.7} is an equality, then $U(Y)\in C^2([0,R))\bigcap C^1([0,R])$, $U_{YY}(R)=\infty$.
\end{rem}
\begin{theorem}\label{T2.4}
Suppose $b>0$, $\alpha>0$, $n<0$, $c\in(0,1)$ and  
\begin{equation}\label{2.10}
\left(1-\frac{\alpha+1}{n}\right)^{\frac{n-1-\alpha}{\alpha}}>\frac{1-c}{\alpha c}.
\end{equation}
Then problem \eqref{1.1}, \eqref{1.2} has a unique classical solution iff
\begin{equation}\label{2.11}
\frac{bR}{2}\leq \left[1-c+c\left(1+Cu^\alpha\zeta_1^\alpha\right)^{\frac{n-1}{\alpha}}\right]\zeta_1 = F(\zeta_1),
\end{equation}
where $\zeta_1$ is the first positive zero of $F'(\zeta)=0$. Moreover, $U(Y)$ satisfies \eqref{2.5}, \eqref{2.6} and
\begin{equation}\label{2.12}
0\leq F^{-1}\left(\frac{bR}{2}\right)\leq \zeta_1 \quad \textrm{for} \quad Y\in [0,R].
\end{equation}
\end{theorem}
\begin{rem}\label{R2.2}
If \eqref{2.11} is a strict inequality, then $U(Y)\in C^2([0,R])$. When \eqref{2.11} is an equality, then $U(Y)\in C^2([0,R))\bigcap C^1([0,R])$, $U_{YY}(R)=\infty$.
\end{rem}
\section{Proofs of the main results}
Integrating $Y.L(U(Y))$ from $0$ to $Y\in(0,R]$, we get the identity
\begin{equation}\label{3.1}
\left [1-c+c \left(1+Cu^\alpha \mid U_Y \mid ^\alpha\right)^{\frac{n-1}{\alpha}}\right]U_{Y}=\frac{bY}{2} \quad \textrm{for} \quad Y\in [0,R].
\end{equation}
Hence $U_Y(Y)>0$ for $Y\in (0,R]$ and \eqref{3.1} is equivalent to
\begin{equation}\label{3.2}
F(U_Y(Y))=\left [1-c+c \left(1+Cu^\alpha U_Y^\alpha\right)^{\frac{n-1}{\alpha}}\right]U_{Y}=\frac{bY}{2} \quad \textrm{for} \quad Y\in [0,R].
\end{equation}

\begin{proof}of Theorem 2.1.
For $c=0$ or $n=1$ we have $F(\zeta)=\zeta$ and \eqref{3.2} becomes
\begin{equation}\label{3.3}
U_Y(Y)=\frac{bY}{2} \quad \textrm{for} \quad Y\in [0,R].
\end{equation}
Integrating \eqref{3.3} from $0$ to $Y\in(0,R]$, we obtain from \eqref{1.2} that $\displaystyle{U(Y)=\frac{b}{4}(Y^2-R^2)}$, which proves Theorem 2.1. 

\end{proof}

\begin {proof} of Theorem 2.2.
Simple computations give us for $\zeta\geq 0$
\begin{equation}\label{3.4}
F'(\zeta)=1-c+c \left(1+Cu^\alpha \zeta^\alpha\right)^{\frac{n-1-\alpha}{\alpha}}\left(1+nCu^\alpha \zeta^\alpha\right),
\end{equation}
\begin{equation}\label{3.5}
F''(\zeta)=c(n-1)Cu^\alpha \zeta^{\alpha-1} \left(1+Cu^\alpha \zeta^\alpha\right)^{\frac{n-1-2\alpha}{\alpha}}\left(\alpha+1+nCu^\alpha \zeta^\alpha\right).
\end{equation}
Thus from \eqref{3.4} we get
\begin{equation}\label{3.6}
F'(\zeta)\geq1-c+c \left(1+Cu^\alpha \zeta^\alpha\right)^{\frac{n-1-\alpha}{\alpha}}>0 \quad \textrm{for} \quad \zeta\geq 0, \quad n\geq0 \quad \textrm{and} \quad c\in (0,1].
\end{equation}
If \eqref{2.4} holds then from \eqref{3.5} we obtain 
\begin{equation}\label{3.7}
F''(\zeta_0)=0 \quad \textrm{for} \quad Cu^\alpha \zeta_0^\alpha=-\frac{\alpha+1}{n}>0,
\end{equation}
\begin{equation}\label{3.8}
F''(\zeta)<0 \quad \textrm{for} \quad  \zeta<\zeta_0 \quad \textrm{i.e. for} \quad Cu^\alpha \zeta^\alpha<-\frac{\alpha+1}{n} \quad \textrm{and}
\end{equation}
\begin{equation}\label{3.9}
F''(\zeta)>0 \quad \textrm{for} \quad  \zeta>\zeta_0 \quad \textrm{i.e. for} \quad Cu^\alpha \zeta^\alpha>-\frac{\alpha+1}{n}.
\end{equation}
From \eqref{3.7}-\eqref{3.9} it follows that the function $F'(\zeta)$ has a global minimum in the interval $[0,\infty)$ at the point $\zeta=\zeta_0$. Since from \eqref{2.4}, we have
\begin{equation}\label{3.10}
F'(\zeta_0)=1-c-\alpha c \left(1-\frac{\alpha+1}{n}\right)^{\frac{n-1-\alpha}{\alpha}},
\end{equation}
\begin{equation}\label{3.11}
F'(\zeta)\geq F'(\zeta_0)>0 \quad \textrm{for} \quad \zeta\geq 0.
\end{equation}
Under the conditions \eqref{2.3}, \eqref{2.4} and from \eqref{3.6}, \eqref{3.1}, it is seen that $F'(\zeta)>0$ and $F(\zeta)$ is strictly monotone increasing. Moreover, 
\begin{equation}\label{3.11a}
F(0)=0, \quad \lim_{\zeta\rightarrow\infty}F(\zeta)=\infty, 
\end{equation}
so that there exits the inverse function $F^{-1}(\zeta)$. From 
\begin{equation}\label{3.12}
(F^{-1}(\zeta))'=\frac{1}{F'(F^{-1}(\zeta))}>0
\end{equation}
it follows that $F^{-1}(\zeta):[0,\infty)\rightarrow[0,\infty)$ is strictly monotone increasing function and \eqref{3.2} is equivalent to
\begin{equation}\label{3.13}
U_Y(Y)=F^{-1}\left(\frac{bY}{2}\right) \quad \textrm{for} \quad Y\in[0,R], 
\end{equation}
where
\begin{equation}\label{3.14}
F^{-1}\left(\frac{bY}{2}\right):[0,R]\rightarrow\left[0,F^{-1}\left(\frac{bR}{2}\right)\right].
\end{equation}
Integrating \eqref{3.13} from $Y\in[0,R)$ to $R$, we get from \eqref{1.2} that the function
\begin{equation}\label{3.15}
U(Y)=-\int^R_Y F^{-1}\left(\frac{bs}{2}\right)ds \quad \textrm{for} \quad Y\in [0,R]
\end{equation}
is the unique classical solution of \eqref{1.1}, \eqref{1.2}, when one of the conditions \eqref{2.3}, \eqref{2.4} is satisfied. Moreover, from the monotonicity of $F^{-1}(\zeta)$ and \eqref{3.13}, estimate \eqref{2.6} holds and Theorem 2.2 is completed.

\end{proof}

\begin{proof} of Theorem 2.3.
From \eqref{2.6a} inequalities \eqref{3.10}, \eqref{3.11} become
\begin{equation}\label{3.16}
F'(\zeta_0)=0, \quad F'(\zeta)>0 \quad \textrm{for} \quad \zeta\geq0, \quad \zeta\neq\zeta_0.
\end{equation}
Thus $F(\zeta)$ is monotone increasing function and the inverse function $F^{-1}(\zeta)$ exists and is monotone increasing, which is evident from \eqref{3.12}. Moreover, \eqref{3.11a} holds and $F(\zeta):[0,\infty)\rightarrow[0,\infty)$, $F^{-1}(\zeta):[0,\infty)\rightarrow[0,\infty)$. Hence, equation \eqref{3.2} is equivalent to \eqref{3.13} and after integration, as in the proof of Theorem 2.2, the unique solution of \eqref{1.1}, \eqref{1.2} is given by \eqref{3.15}.

Let us analyze the regularity of $U(Y)$. It is clear that $U(Y)\in C^1([0,R])$ and from the monotonicity of $F^{-1}$ the estimate \eqref{2.6} holds. Since \eqref{2.7} is equivalent to
\begin{equation}\label{3.17}
F(\zeta_0)\geq \frac{bR}{2}, \quad \textrm{i.e.,} \quad \zeta_0\geq F^{-1}\left(\frac{bR}{2}\right)
\end{equation}
and after differentiating \eqref{3.13}, we get from \eqref{3.12}, \eqref{3.8} the estimates 
\begin{equation}\label{3.18}
U_{YY}(Y)=\frac{\partial }{\partial Y}\left(F^{-1}\left(\frac{bY}{2}\right)\right)=\frac{b}{2F'(F^{-1}\left(\frac{bY}{2}\right))}\leq\frac{b}{2F'(F^{-1}\left(\frac{bR}{2}\right))} 
\end{equation}
\[
\leq \frac{b}{2F'(F^{-1}\left(\zeta_0\right))}=\infty \quad \textrm{for} \quad Y\in [0,R].
\]
Thus
\begin{equation}\label{3.19}
0< U_{YY}(Y)<\infty \quad \textrm{for} \quad Y\in [0,R) \quad \textrm{and} \quad U_{YY}(R) = \infty,
\end{equation}
when \eqref{2.7}, equivalently \eqref{3.17}, is a strict inequality, because 
\[\displaystyle{\frac{b}{2F'(F^{-1}\left(\frac{bR}{2}\right))}< \frac{b}{2F'(F^{-1}\left(\zeta_0\right))}=\infty }.\] 
However, when \eqref{2.7}, equivalently \eqref{3.17}, is equality, then \[\displaystyle{U_{YY}(R)=\frac{b}{2F'(F^{-1}\left(\frac{bR}{2}\right))}=\frac{b}{2F'(F^{-1}\left(\zeta_0\right))}=\infty }.\]
Remark 2.1. and the proof of \eqref{2.7} is completed.

Suppose that \eqref{2.8} holds. As in the proof of \eqref{2.7} $U(Y)\in C^1([0,R])$ is unique solution of \eqref{1.1}, \eqref{1.2} satisfying \eqref{2.5}, \eqref{2.6}. Conditions \eqref{2.8}, \eqref{2.9} are equivalent to 
\begin{equation}\label{3.20}
F(\zeta_0)< \frac{bR}{2}, \quad F(\zeta_0)= \frac{bY_0}{2}, \quad \textrm{i.e.,} \quad \zeta_0< F^{-1}\left(\frac{bR}{2}\right), \quad \zeta_0= F^{-1}\left(\frac{bY_0}{2}\right)
\end{equation}
and therefore $Y_0\in (0,R)$.\\
From the monotonicity of $F^{-1}$ we get from \eqref{3.20}
\begin{equation}\nonumber
F^{-1}\left(\frac{bY}{2}\right)< F^{-1}\left(\frac{bY_0}{2}\right)=\zeta_0 \quad \textrm{for} \quad Y\in [0,Y_0),
\end{equation}
\begin{equation}\label{3.21}
 F^{-1}\left(\frac{bY}{2}\right)> F^{-1}\left(\frac{bY_0}{2}\right)=\zeta_0 \quad \textrm{for} \quad Y\in (Y_0,R].
\end{equation}
Repeating the estimate \eqref{3.18}, from \eqref{3.21} and the monotonicity of $F'(\zeta)$, see \eqref{3.8}, \eqref{3.9}, we have
\begin{equation}\nonumber
U_{YY}(Y_0)=\infty,\quad 0< U_{YY}(Y) < U_{YY}(Y_0)=\infty \quad \textrm{for} \quad Y\in [0,R]\setminus {Y_0},
\end{equation}
i.e., $U(Y)\in C^2([0,R] \setminus {Y_0})\bigcap C^1([0,R])$, which proves \eqref{2.8} and Theorem 2.3.

\end{proof}

\begin{proof} of Theorem 2.4.
Under conditions \eqref{2.10}, \eqref{3.10} it follows that $F'(\zeta)$ has a strictly negative minimum in the interval $[0,\infty)$ at the point $\zeta_0$. Since
\begin{equation}\label{3.22}
F'(0)=0, \quad \lim_{\zeta\rightarrow\infty}F'(\zeta)=1-c>0 
\end{equation}
the function $F'(\zeta)$ has two positive roots $0<\zeta_1<\zeta_2<\infty$, $F'(\zeta_i)=0, \quad i=1,2$, $\zeta_0\in(\zeta_1,\zeta_2)$. Moreover, from \eqref{3.22} and \eqref{3.8} it follows that
\begin{equation}\label{3.23}
F'(\zeta)>0, \quad F''(\zeta)<0 \quad \textrm{for} \quad \zeta\in [0,\zeta_1), \quad F'(\zeta_1)=0,
\end{equation}
i.e., $F'(\zeta)$ is strictly monotone decreasing, while $F(\zeta)$ is strictly monotone increasing for 
$\zeta\in [0,\zeta_1)$. Hence, the inverse function 
\begin{equation}\label{3.24}
F^{-1}(\zeta):[0,\zeta_1)\rightarrow[0,F^{-1}(\zeta_1)]
\end{equation}
is well defined and from \eqref{3.12}, \eqref{3.23} is strictly monotone increasing.

\underline{\textit{Sufficiency:}} From \eqref{2.11} we get 
\begin{equation}\label{3.25}
\zeta_1\geq F^{-1}\left(\frac{bR}{2}\right) \quad \textrm{and} \quad F^{-1}\left(\frac{bY}{2}\right):[0,R]\rightarrow \left[0,F^{-1}\left(\frac{bR}{2}\right)\right].
\end{equation}
Since $\displaystyle{\left[0,F^{-1}\left(\frac{bR}{2}\right)\right]\subset [0,\zeta_1)}$, from \eqref{3.24} equation \eqref{3.2} is equivalent to \eqref{3.13} and after integration $U(Y)\in C^1([0,R)$ satisfies \eqref{3.15}. Moreover, from the monotonicity of $F^{-1}(\zeta)$ and \eqref{3.25}, the estimate \eqref{2.12} is satisfied because 
\[
0\leq U_Y(Y) \leq U_Y(R)=F^{-1}\left(\frac{bR}{2}\right)\leq \zeta_1.
\] 
Repeating the estimate \eqref{3.18} we get 
\begin{equation}\nonumber
0<U_{YY}(Y)<\infty \quad \textrm{for} \quad Y\in [0,R],
\end{equation}
when \eqref{2.11}, equivalently \eqref{3.25}, is a strict inequality, i.e., $U(Y)\in C^2([0,R])$.

If \eqref{2.11}, equivalently \eqref{3.25}, is an equality, then 
\begin{equation}\nonumber
U_{YY}(R)=\frac{b}{2F'(F^{-1}\left(\frac{bR}{2}\right))}= \frac{b}{2F'(\zeta_1)}=\infty, \quad  U(Y)\in C^2([0,R))\cap C^1([0,R]),
\end{equation}
which proves the sufficiency in Theorem 2.4 and Remark 2.2.

\underline{\textit{Necessity:}} Suppose $U(Y)\in C^2((0,R))\cap C^1([0,R])$ is a classical solution of \eqref{1.1}, \eqref{1.2}. Then $U(Y)$ satisfies \eqref{3.2}. We assume by contradiction that \eqref{2.11} fails, i.e., 
\begin{equation}\label{3.26}
\frac{bR}{2}> F(\zeta_1).
\end{equation}
If $\displaystyle{Y_0=\frac{2F(\zeta_1)}{b}}$, then from \eqref{3.26}, we have $Y_0\in (0,R)$ and from \eqref{3.23}, \eqref{3.24}
\begin{equation}\label{3.27}
\zeta_1=F^{-1}\left(\frac{bY_0}{2}\right).
\end{equation}
Thus \eqref{3.27} and \eqref{3.18} for $Y=Y_0$ give us
\begin{equation}\label{3.28}
U_{YY}(Y_0)=\frac{b}{2F'(F^{-1}\left(\frac{bY_0}{2}\right))}= \frac{b}{2F'(\zeta_1)}=\infty.
\end{equation}
Since $Y_0$ is an interior point of the interval $(0,R)$, it follows from \eqref{3.28} that $U(Y)$ is not a classical $C^2(0,R)$ smooth solution of \eqref{1.1}, \eqref{1.2} and Theorem 2.4 is completed.

\end{proof}
\section{Applications}
\subsection{Case c=1} In the applications of polymer pipe flows, the special case of $c=1$ is important ($\mu_{\infty}\rightarrow 0$ \cite{Bird:2002}).
 \begin{theorem}\label{T4.1}
 Suppose $b>0$, $\alpha>0$, $n\in \R$ and $c=1$. Then
 \begin{itemize}
 \item{ (a)} for
 \begin{equation}\label{4.1}
\quad n>0
 \end{equation}
 problem \eqref{1.1}, \eqref{1.2} has a unique classical solution $U(Y)\in C^2([0,R])$ satisfying \eqref{2.5}, \eqref{2.6};
 \item{ (b)} for
 \begin{equation}\label{4.2}
  \quad n=0
  \end{equation}
 problem \eqref{1.1}, \eqref{1.2} has a unique classical solution $U(Y)\in C^2((0,R))\cap C^1([0,R])$ iff
   \begin{equation}\label{4.3}
 \frac{bR}{2}\leq Cu^{-1} 
 \end{equation}
and $U(Y)$ satisfies \eqref{2.5}, \eqref{2.6};
\item{ (c)} for
  \begin{equation}\label{4.4}
 \quad n<0
  \end{equation}
 problem \eqref{1.1}, \eqref{1.2} has a unique classical solution $U(Y)\in C^2((0,R))\cap C^1([0,R])$ iff
   \begin{equation}\label{4.5}
 \frac{bR}{2}\leq Cu^{-1} \left(\frac{n-1}{n}\right)^{\frac{n-1}{\alpha}}\left(-\frac{1}{n}\right)^{\frac{1}{\alpha}}=F\left(Cu^{-1}\left(-\frac{1}{n}\right)^{\frac{1}{\alpha}}\right)
 \end{equation}
and $U(Y)$ satisfies \eqref{2.5}, \eqref{2.6}.
\end{itemize}
 \end{theorem}
 \begin{rem}\label{R4.1}
 If \eqref{4.3} or \eqref{4.5} is strict inequality then $U(Y)\in C^2([0,R])$. When \eqref{4.3} or \eqref{4.5} is equality then $U(Y)\in C^2([0,R))\cap C^1([0,R])$ and $U_{YY}(R)=\infty$.
 \end{rem}

\begin{proof}
\begin{itemize}
\item {(a)} 
Suppose \eqref{4.1} holds. Then from \eqref{3.6}
\begin{equation}\label{4.6}
F'(\zeta)\geq \left(1+Cu^\alpha\zeta^\alpha \right)^{\frac{n-1-\alpha}{\alpha}}>0 \quad \textrm{for} \quad \zeta\geq 0
  \end{equation}
  and $F(\zeta)$ is strictly monotone increasing function. 
  
 Since $F(0)=0$, $\lim_{\zeta\rightarrow\infty}F(\zeta)=\infty$, there exists the inverse function $F^{-1}(\zeta):[0,\infty)\rightarrow[0,\infty)$, which is strictly monotone increasing function from \eqref{3.12} and \eqref{4.6}. Thus problem \eqref{1.1}, \eqref{1.2} is equivalent to \eqref{3.2} and after integration, $U(Y)$ satisfies \eqref{2.5}, \eqref{2.6}.
\item  {(b) }
Suppose \eqref{4.2} holds. Here, again \eqref{4.6} is valid, from which $F(\zeta)$ is strictly monotone increasing function for $\zeta\geq 0$. Since $F(0)=0$, $\lim_{\zeta\rightarrow\infty}F(\zeta)=Cu^{-1}=\sup_{\zeta\geq 0} F(\zeta)$, there exists the inverse function
 \begin{equation}\label{4.7}
F^{-1}(\zeta):[0,Cu^{-1})\rightarrow[0,\infty),
  \end{equation}
  which is strictly monotone increasing function from \eqref{3.12} and \eqref{4.6}.
  
  \underline{\textit{Sufficiency:}} From \eqref{4.3} 
   $\displaystyle{F^{-1}\left(\frac{bY}{2}\right):\left[0,\frac{bR}{2}\right]\rightarrow \left[0,F^{-1}\left(\frac{bR}{2}\right)\right]}$ is well defined for  $Y\in[0,R]$ because $\displaystyle{ \left[0,\frac{bR}{2}\right]\subset [0,Cu^{-1}]}$. The rest of the proof is identical with the proof of \textit{(a)} of Theorem 4.1.
   
\underline{\textit{Necessity:}} If $U(Y)\in C^2((0,R))\cap C^1([0,R])$ is a classical solution of \eqref{1.1}, \eqref{1.2}, we suppose by contradiction that $\displaystyle{\frac{bR}{2}> Cu^{-1}}$. 

Then $\displaystyle{Y_0=\frac{2}{bCu}\in (0,R)}$ and from \eqref{3.2} and the monotonicity of $F(\zeta)$ we have for $Y\in (Y_0,R)$ the following impossible chain of inequalities
\begin{equation}\nonumber
Cu^{-1}=\frac{bY_0}{2}<\frac{bY}{2}=F(U_Y(Y))\leq\sup_{\zeta\geq 0}F(\zeta)=Cu^{-1}.
\end{equation}
\item{(c)}
Suppose \eqref{4.4} holds. From \eqref{3.4}, \eqref{3.5} we get
    \begin{eqnarray}
    F'(\zeta)=\left(1+Cu^\alpha \zeta^\alpha\right)^{\frac{n-1-\alpha}{\alpha}}\left(1+nCu^\alpha \zeta^\alpha\right),\\
 F''(\zeta) =(n-1)Cu^\alpha \zeta^{\alpha-1} \left(1+Cu^\alpha \zeta^\alpha\right)^{\frac{n-1-\alpha}{\alpha}}\left(\alpha+1+nCu^\alpha \zeta^\alpha\right).   
    \end{eqnarray}
 and $F'(\zeta)> 0$ for $\displaystyle{Cu^\alpha \zeta^\alpha < -\frac{1}{n}}$, $F'(\zeta_1)=0$ for $\displaystyle{Cu^\alpha \zeta^\alpha = -\frac{1}{n}}$
 and $F'(\zeta)< 0$ for $\displaystyle{Cu^\alpha \zeta^\alpha > -\frac{1}{n}}$. Thus the function $F(\zeta)$ has a global maximum in the interval $[0,\infty)$ at the point $\displaystyle{\zeta_1=Cu^{-1}\left(-\frac{1}{n}\right)^{\frac{1}{\alpha}}}$ and 
\[F(\zeta_1)=Cu^{-1}\left(\frac{n-1}{n}\right)^{\frac{n-1}{\alpha}}\left(-\frac{1}{n}\right)^{\frac{1}{\alpha}}> 0.
\]
 
 The rest of the proof is identical with the proof of Theorem 2.4 and we omit it. The proof of remark 4.1 is identical with the proof of Remarks 2.1 and 2.2.
 \end{itemize}
 
\end{proof}

In practice the necessary and sufficient condition \eqref{4.5} corresponds to three different scenarios for solution existence of \eqref{1.1}, \eqref{1.2}, which depend on the values of the parameters $n$, $\alpha$, $Cu$, $b$ and $R$. Here, we present some examples to illustrate each of these three cases.
\begin{itemize}
 \item{ (i)} if \eqref{4.5} is strict inequality, then $\zeta_1$ is not reached for any $Y\in[0,R]$, i.e., $\zeta(R) < \zeta_1$ and $\zeta_1=\zeta(Y_1)$, where $Y_1 > R$. For example, at $n=-10$, $\alpha=10$, $Cu=b=R=1$, we have $\displaystyle{Cu^{-1} \left(\frac{n-1}{n}\right)^{\frac{n-1}{\alpha}}\left(-\frac{1}{n}\right)^{\frac{1}{\alpha}}=0.7153}$, $\displaystyle{\frac{bR}{2}=0.5}$ and \eqref{4.5} is fulfilled. Moreover, $\zeta_1=0.7943$ and $Y_1=1.4305 > R=1$. This means that if $\zeta \leq \zeta_1$, then the classical solution exists everywhere inside the pipe.
From the equation $\displaystyle{F(\zeta)=\frac{bY}{2}}$, the inverse function $\displaystyle{F^{-1}\left(\frac{bY}{2}\right)}$ is constructed numerically using the computational algebra package MAPLE. Two different functions $\displaystyle{F^{-1}\left(\frac{bY}{2}\right )}$ are found:
\[ F^{-1}_1\approx 0.02232Y^9-0.07589Y^7+0.08834Y^5-0.03763Y^3+0.5044Y,
\]
 \[
F^{-1}_2\approx \frac{16.03 Y}{(Y^2)^{2.064}+16.73 (Y^2)^{0.5535}}
\]
The last function is not bounded at $Y=0$, i.e.,  $F^{-1}_2(0)=U_{2Y}(0)\rightarrow\infty$, which means that $U_{2Y}\neq 0$ does not fulfill the condition \eqref{1.2}.
Then, only the function $F^{-1}_1$ corresponds to the velocity $U(Y)$:
\[
U_1(Y)=0.00223Y^{10}-0.00949Y^8+0.01472Y^6-0.0094Y^4+0.2522Y^2-0.25026,
\]
which is solution of \eqref{1.1}, \eqref{1.2}.
\item{(ii)} if \eqref{4.5} is equality, then $\displaystyle{F(\zeta_1)=\frac{bR}{2}}$ and $\displaystyle{\zeta_1=F^{-1}\left(\frac{bR}{2}\right)}$ is reached at $Y_1=R$. For example, at $n=-5$, $\alpha=3.9$, $Cu=b=R=1$, we have $\zeta_1=0.662$, $Y_1=R=1$ and $\displaystyle{F(\zeta_1)=\frac{bR}{2}}=0.5$. The inverse function 
\[
F^{-1}\left(\frac{bY}{2}\right)\approx 
0.7761Y^9-1.301Y^7+0.798Y^5-0.1547Y^3+0.5085Y
\]
is again constructed numerically and 
\[
U(Y)=0.07761Y^{10}-0.16262Y^8+0.133Y^6-0.03868Y^4+0.25425Y^2-0.26356
\]
 which is solution of \eqref{1.1} and \eqref{1.2}.
\item{(iii)} if \eqref{4.5} is not fulfilled, then $\zeta_1$ is reached at some inner point $Y_1\in[0,R)$, such that $\zeta_1=\zeta(Y_1)<\zeta(R)$. This means that the equation \eqref{1.1} together with the first condition in \eqref{1.2} has solution, which is not a solution of the problem \eqref{1.1}, \eqref{1.2}. For example, at $n=-3$, $\alpha=2$, $Cu=b=R=1$, it is obtained that $\zeta_1=0.325$, $Y_1=0.65 < R=1$ and $\displaystyle{F(\zeta_1)=0.325 < \frac{bR}{2}}=0.5$. 
Thus if $\zeta \leq \zeta_1$, the classical solution exists only in $[0, Y_1]$. There is no solution in the rest of the pipe cross section, i.e., in $(Y_1, R]$. In $[0, Y_1]$ the inverse function is
\[
F^{-1}\left(\frac{bY}{2}\right )\approx 48.1Y^9-33.83Y^7+8.55Y^5-0.4777Y^3+0.5168Y.
\] 
This function corresponds to a single velocity profile \[U(Y)= 4.81Y^{10}-4.2287Y^8+1.425Y^6-0.11942Y^4+0.2584Y^2+const.,\]
with an unknown constant depending on an eventual additional condition at $Y=Y_1$ or at $Y=0$.
\end{itemize}
The given examples are arbitrary chosen with the only aim to illustrate the importance of the condition \eqref{4.5}.
\subsection{Revision of unsteady case} Another application of Theorem 2.4 and Theorem 4.1. is connected with the unsteady Carreau-Yasuda flow in a pipe. We revise the results in \cite{Kutev:2021c}, where the following problem is studied
\begin{equation}\label{4.8}
8\beta^2W_T - \frac{1}{Y}\frac{\partial}{\partial Y}\left \{ \left [1-c+c \left(1+Cu^\alpha \mid W_Y \mid  ^\alpha\right)^{\frac{n-1}{\alpha}}\right] Y W_{Y} \right \}=f(T),
\end{equation}
\begin{equation}\nonumber
\textrm{for} \quad T>0, \quad Y\in(0, R),\qquad 
\end{equation}  
\begin{equation}\label{4.9}
W_Y(T,0) = W(T,R) = 0 \quad \textrm{for} \quad T \geq 0,  
\quad
W(0,Y)=\Psi(Y) \quad \textrm{for} \quad  Y\in[0, R],
\end{equation}
with  $\Psi(Y)\in C^4([0,R])$ satisfying the compatibility conditions:
\begin{equation}\label{4.10}
\Psi'(0) = \Psi'''(0)=\Psi(R) = 0. 
 \end{equation}
and $f \in C^1([0,\infty))$ bounded as
 \begin{equation}\label{4.11}
\sup_{T \ge 0}\mid f(T) \mid < \infty.
 \end{equation}
 The problem parameters are: $c\in[0,1)$, $\alpha>0$, $n\in \R$, $\beta>0$. 
 
As it is proved in \cite{Kutev:2021c}, under the condition
\begin{equation}\label{4.12}
\alpha\left(1-\frac{\alpha+1}{n}\right)^{\frac{n-1-\alpha}{\alpha}} > \frac{1-c}{c}, \quad n<0
\end{equation}
equation \eqref{4.8} becomes forward-backward parabolic equation. If \eqref{4.12} holds, then the function $h(\eta)$
\begin{equation}\label{4.13}
h(\eta)=1-c+c\Phi(\eta), \quad \textrm{where}
\end{equation}
\begin{equation}\label{4.14}
\eta=\left | \frac{\partial W} {\partial Y}\right | ,  \quad \Phi(\eta) = (1-n)\left(1+Cu^\alpha\eta^\alpha\right)^{\frac{n-1-\alpha}{\alpha}} + n\left(1+Cu^\alpha\eta^\alpha\right)^{\frac{n-1}{\alpha}}
\end{equation}
has two positive roots, $0<\eta_1<\eta_2$, such that \eqref{4.8} is forward parabolic equation for $\eta \in [0,\infty) \setminus (\eta_1,\eta_2)$ and backward parabolic one for  $\eta \in (\eta_1,\eta_2)$.

We recall Theorem 3 in \cite{Kutev:2021c}, where sufficient conditions for existence of global solutions of \eqref{4.8}-\eqref{4.11} are given.
\begin{theorem}\label{T4.2} (Theorem 3 in \cite{Kutev:2021c})
Suppose $W(T,Y)\in C^3\left([0,T_m] \times [0,R] \right)$ is a solution of \eqref{4.8} - \eqref{4.11}, $0\leq T_m\leq \infty$, $\alpha>0$, $\beta>0$, $n < 0$, $c\in (0,1)$ and \eqref{4.12} hold. If $K_1<\eta_1$, where $K_1$ is given in \eqref{4.15} and $\eta_1$ is the first positive root of $h(\eta)$, then $T_m=\infty$ and the estimates
\begin{equation}\nonumber
\left|W(T,Y)\right|\leq K_1\left(R^2-Y^2\right),
\end{equation}
\begin{equation}\nonumber
\left|W_{Y}(T,Y)\right| \leq K_1 
\end{equation}
are satisfied for $T \geq 0, Y \in [0, R]$. The constant $K_1$ is defined in the following way
\begin{equation}\label{4.15}
K_1 = R M\left [ 1-c+c\left(1+Cu^\alpha\lim_{m\rightarrow\infty}(C_m)^\alpha\right)^{\frac{n-1}{\alpha}}\right ]^{-1}, \quad \textrm{where}
\end{equation}
\begin{equation}\nonumber
M=\max \Bigg\lbrace \sup_{T \ge 0} \mid f(T) \mid, \frac{1}{R}\sup_{Y \in [0, R]}\mid \Psi'(Y)\mid, 2\sup_{Y \in [0, R]} \frac{\mid \Psi(Y)\mid}{R^2-Y^2}\Bigg\rbrace,
\end{equation}
\begin{equation}\nonumber
C_m=R M\left [ 1-c+c\left(1+Cu^\alpha(C_{m-1})^\alpha\right)^{\frac{n-1}{\alpha}}\right ]^{-1}, \quad m=1,2,\cdots
\end{equation}
\[
C_0=\frac{R M}{1-c} 
\]
\end{theorem}
By means of Theorem 2.4 and Theorem 4.1, we have the following precise variant of the upper Theorem 4.2. For this purpose we take the right-hand side of \eqref{1.1} as $b=\sup_{T\geq 0}\left|f(T)\right |$, $\alpha>0$, $\beta>0$, $n < 0$, $c\in (0,1)$ and when \eqref{2.10}, \eqref{2.11} hold.
\begin{theorem}\label{T4.3} 
Suppose $\alpha>0$, $\beta>0$, $n < 0$, $c\in (0,1]$ and \eqref{2.10}, \eqref{2.11} hold. If
\begin{equation}\label{4.16}
-U(Y)\leq \Psi(Y) \leq U(Y) \quad \textrm{for} \quad  Y\in[0, R],
\end{equation}
\begin{equation}\label{4.17}
\max \Bigg\lbrace F^{-1}\left( \frac{R}{2} \sup_{T \ge 0} \mid f(T) \mid \right), \sup_{Y \in [0, R]}\mid \Psi'(Y)\mid \Bigg\rbrace < \zeta_1 
\end{equation}
then the problem \eqref{4.8} - \eqref{4.11} has a unique global classical solution $W(T,Y) \in C^3\left([0,\infty) \times [0,R] \right)$ and the estimate
\begin{equation}\label{4.18}
\left|W_Y(T,Y)\right| \leq \zeta_1 \quad \textrm{for} \quad T\geq 0, \quad Y\in[0, R]
\end{equation}  
is satisfied, where $\zeta_1$ is well defined in Theorem 2.4.
\begin{proof}
From \eqref{4.6} and the comparison principle, Theorem 1 in \cite{Kutev:2021c}, we get the estimates
\begin{equation}\label{4.19}
-U(Y)\leq W(T,Y) \leq U(Y) \quad \textrm{for} \quad T\geq 0, \quad  Y\in[0, R],
\end{equation}
\begin{equation}\label{4.20}
\left|W_Y(T,Y)\right| \leq \left|U_Y(Y)\right| =  \left|F^{-1}\left( \frac{R}{2} \sup_{T \ge 0} \mid f(T) \mid \right)\right| \quad \textrm{for} \quad T\geq 0
\end{equation}
The rest of the proof follows from \eqref{4.17}, \eqref{4.20} and the global gradient estimates, Theorem 2 in \cite{Kutev:2021c}.

\end{proof}
\end{theorem}

\subsection*{A Appendix} 
We will prove that $\zeta_1\equiv \eta_1$, where $\zeta_1$ is the first positive zero of $F'(\zeta)=0$ and $\eta_1$ is the first positive zero of $h(\eta)=0$.

From \eqref{3.4} we define $Z=1+Cu^\alpha\zeta_1^\alpha$ as a solution of the equation
\begin{equation}\nonumber
1+n Cu^\alpha\zeta_1^\alpha=-\left(\frac{1-c}{c}\right)\left(1+Cu^\alpha\zeta_1^\alpha\right)^{\frac{1+\alpha-n}{\alpha}}, \quad \textrm{i.e.,}
\end{equation}
\begin{equation}\nonumber
\left(\frac{1-c}{c}\right)Z^{\frac{1+\alpha-n}{\alpha}}+n Z-n+1 =0.
\end{equation}
In the same way, $\eta_1$ is a solution of the equation 
\begin{equation}\nonumber
1+n Cu^\alpha\eta_1^\alpha=-\left(\frac{1-c}{c}\right)\left(1+Cu^\alpha\eta_1^\alpha\right)^{\frac{1+\alpha-n}{\alpha}}, \quad \textrm{i.e.,}
\end{equation}
and for $X=1+Cu^\alpha\eta_1^\alpha$, we get
\begin{equation}\nonumber
\left(\frac{1-c}{c}\right)X^{\frac{1+\alpha-n}{\alpha}}+n X-n+1 =0.
\end{equation}
Since $\zeta_1<\zeta_0$, where $\displaystyle{Cu^\alpha\zeta_0^\alpha=-\frac{\alpha+1}{n}}$ and $\eta_1<\eta_0$, where $\displaystyle{Cu^\alpha\eta_0^\alpha=-\frac{\alpha+1}{n}}$, we conclude that $\zeta_1=\eta_1$. Here $\eta_0$ is the maximum point for the function $h(\eta)$, i.e.,
 \[
h'(\eta_0)=0, h'(\eta)=c(n-1)Cu^\alpha \left(1+Cu^\alpha\eta^\alpha\right)^{\frac{n-1-2\alpha}{\alpha}} \left(1+\alpha+nCu^\alpha\eta^\alpha\right)\eta^{\alpha-1}.
\]

\subsection*{Acknowledgment}

N.K. has been partially supported by the National Scientific Program "Information and Communication Technologies for a Single Digital Market in Science, Education and Security (ICTinSES)", contract No D01-–205/23.11.2018, financed by the Ministry of Education and Science in Bulgaria and by the Grant
No BG05M2OP001-1.001-0003-C01, financed by the Science and Education for Smart Growth
Operational Program (2018-2023).

\end{document}